\documentclass[11pt]{article}
\usepackage{array, amsmath, amssymb}

\begin{document}
\begin{center}
{\Large\bf A further generalization of the Emden--Fowler equation}\\
\vspace{0.5cm}
{\bf Shinji Tanimoto}\\ 
\vspace{0.5cm}
Department of Mathematics,
University of Kochi,\\
Kochi 780-8515, Japan\footnote{Former affiliation}.  \\
\end{center}
\begin{abstract}
A generalization of the Emden--Fowler equation is presented and its solutions are investigated. 
This paper is devoted to asymptotic behavior of its solutions.  The procedure is entirely based on a previous paper by the author. 
\end{abstract}
\vspace{0.7cm}
%%%%   1  %%%%%%%%%%%%%%%%%%%%%%%%%%%%%%%%%
{\large \bf 1. Introduction} \\  
\\
\indent
The Emden--Fowler equation is a second-order differential equation taking the form
\begin{eqnarray}
y'' = \phi(t) |y|^{\lambda}\,{\rm sgn }\,y,
\end{eqnarray}
where $\phi$ is a continuous function and $\lambda$ is a positive number. This type of equations plays an
important role in many areas of theoretical physics. As a particular case it includes the Thomas--Fermi equation
in atomic physics
\begin{eqnarray*}
y'' = t^{-1/2}  y^{3/2}.
\end{eqnarray*}
Eq.(1) is linear when $\lambda = 1$, and it is superlinear or sublinear, when $\lambda > 1$ or $\lambda < 1$,
respectively. \\
\indent
Many authors have mainly investigated its asymptotic solutions $y(t)$ as $t \to \infty$.
We refer only to [2] and the references therein for such results,
because our present results do not share so much relationship to them. \\
\indent
In this paper we propose a more general equation and study its asymptotic solutions. It is
an $n$th--order differential equation of the form
\begin{eqnarray}
y^{(n)} = F(t, y),
\end{eqnarray}
where $F(t,y)$ is a continuous function of two variables $t$ and $y$. This equation is also a generalization of the equation
\begin{eqnarray}
y^{(n)}= \phi(t) |y|^{\lambda}\,{\rm sgn }\,y,
\end{eqnarray}
considered in [1]; i.e., $F(t,y)= \phi(t) |y|^{\lambda}\,{\rm sgn }\,y$. \\
\indent
Our aim is to study asymptotic behavior of its solutions $y(t)$ of Eq.(2) as $t \to \infty$ under assumptions imposed on 
the function $F$. The procedure is based on [1], where asymptotic properties of functions are discussed. As for the results of [1], 
we briefly recall them together with the basic definitions in the next section. 
Although the notation here is different from that of [1], the substance is all the same. \\
\\
%%%%   2  %%%%%%%%%%%%%%%%%%%%%%%%%%%%%%%%%
{\large \bf 2. Asymptotic behavior of functions} \\  
\\
\indent 
We assume that all the functions are continuous and real--valued, whose domains are intervals of the type $[a, \infty)$
depending on functions. According to aymptotic properties of functions, we divide such functions into 
three categories; $\mathcal S_1$, $\mathcal S_2$ and $\mathcal S_3$. 
In the following '$\limsup$' and '$\liminf$' mean $\limsup_{t \to +\infty}$ and $\liminf_{t \to +\infty}$, respectively.
\begin{itemize}
\item[($\mathcal S_1$)] $\mathcal S_1$ denotes the set of all functions $f$ such that
\begin{eqnarray*} \liminf \frac{|f(t)|}{t^p} = \infty,
\end{eqnarray*} for all real numbers $p$. Some typical examples include $e^{\sqrt t}$, $e^{t^2}$.     
It is obvious that in case of $f \in \mathcal S_1$ we have $|f|^{\nu} \in \mathcal S_1$ for $\nu > 0$.
\item[($\mathcal S_2$)] $\mathcal S_2$ denotes the set of all functions $f$ such that
\begin{eqnarray*} \limsup \frac{|f(t)|}{t^p} = 0,
\end{eqnarray*} for all real numbers $p$. Some typical examples include $e^{-\sqrt t}$, $e^{-t^2}$.             
\item[($\mathcal S_3$)] $\mathcal S_3$ denotes the set of all functions $f$, not belonging to $\mathcal S_1$ nor
to $\mathcal S_2$.
\end{itemize}
With each function $f \in \mathcal S_3$ we associate two values $\Pi(f)$ and $\Xi(f)$. The value
$\Pi(f)$ is characterized by the number such that for all $\epsilon > 0$
\begin{eqnarray*}\liminf |f(t)|/t^{\Pi(f) -\epsilon} = \infty ~~{\rm and} ~~ \liminf |f(t)|/t^{\Pi(f) + \epsilon} = 0.
\end{eqnarray*}
If there exists no such a number (i.e., $\liminf |f(t)|/t^p = 0$ for all $p$), put $\Pi(f) = - \infty$. \\
\indent
On the other hand, the value $\Xi(f)$ is characterized by the number such that for all $\epsilon > 0$
\begin{eqnarray*}\limsup |f(t)|/t^{\Xi(f) -\epsilon} = \infty ~~{\rm and}~~ \limsup |f(t)|/t^{\Xi(f) + \epsilon} = 0.
\end{eqnarray*}
If there exists no such a number (i.e., $\limsup |f(t)|/t^p = \infty$ for all $p$), put $\Xi(f) = + \infty$.\\
It is obvious that $\Xi(f) \ge \Pi(f)$ for every $f \in {\mathcal S_3}$. As its example, taking 
$f(t) = t + e^t \cos t$, we have $\Pi(f) = 1$ and $\Xi(f) = + \infty$.  For a polynomial $f$ of degree $d$, we have
$\Pi(f) = \Xi(f) = d$.  \\
\indent
Note that when $\Pi(f) < 0$ or $f \in \mathcal S_2$ there exists an increasing sequence $\{t_k\}$ such that $t_k \to \infty$ and
$\lim_{k \to \infty} f(t_k) = 0$. Furthermore, the smaller $\Pi(f)$ is, the faster $f(t_k)$ approaches to zero.\\
\indent
The following asymptotic behavior concerning the derivatives will be useful, when we consider some types of
differential equations. We implicitly assume the existence of solutions with initial conditions, for which 
we discuss their asymptotic behavior.\\
\\
{\bf Theorem 1.} (see [1, Th3]) {\it Let $f$ be a continuously differentiable function.
Then the following properties hold.}
\begin{itemize}
\item[{\rm (i)}] {\it If both $f$ and its derivative $f'$ belong to {$\mathcal S_3$}, then} $\Pi(f') \le \Pi(f) -1$.
\item[{\rm (ii)}] {\it If $f$ belongs to {$\mathcal S_2$}, then either $f' \in {\mathcal S_2}$ or 
$f' \in {\mathcal S_3}$ with $\Pi(f') = - \infty$}.
\item[{\rm (iii)}]  {\it If the derivative $f'$ belongs to {$\mathcal S_1$}, then} $f \in {\mathcal S_1}$.\\
\end{itemize}
In particular note that (i) also implies $\Pi(f') = -\infty$ for functions with $\Pi(f) = -\infty$. The proofs are straightforward and can be found in [1].  \\
\\
\\
%%%%   3  %%%%%%%%%%%%%%%%%%%%%%%%%%%%%%%%%
{\large \bf 3. Asymptotic behavior of solutions}  
\\
\\
\indent
We study a particular case of Eq.(2) that corresponds to the superlinear case of Eq.(1), i.e., $\lambda >1$. 
The following assumptions on the function $F$ reflect such a case. Indeed we can easily see that the function 
\begin{eqnarray*}
F(t,y)= \phi(t) |y|^{\lambda}\,{\rm sgn }\,y ~~ {\rm for}~ \lambda >1
\end{eqnarray*}
fulfills these assumptions for appropriate functions $\phi$ and appropriate numbers $\nu$ ($1 < \nu <\lambda$), 
by virtue of [1, Th4 and Th5]. We suppress the variable
$t$ for expressing a function (of $t$) itself such as $F(\cdot, y(\cdot))$, for example.  \\
\\
({\bf Assumptions})  
\begin{itemize}
\item[{\rm (i)}]  For some number $\nu > 1$ the function $F$ satisfies
\begin{eqnarray*}
\frac{F(\cdot, f(\cdot))}{(f(\cdot))^{\nu}} \in \mathcal S_1~{\rm for~all}~ f \in \mathcal S_1.
\end{eqnarray*}
\item[{\rm (ii)}]  For all $f \in \mathcal S_3$, $F(\cdot, f(\cdot)) \in \mathcal S_3$ and 
\begin{eqnarray*}
\Pi(F(\cdot, f(\cdot))) \ge \mu \Pi(f) + c
\end{eqnarray*}
holds, where $\mu$ and $c$ are two constants such that $\mu > 1$ and $c > -n$, $n$ being the order of Eq.(2). \\
\end{itemize}
As for the right-hand side of Eq.(3) it is shown in [1] that, when $\Pi(\phi) = \Xi(\phi)$ holds and hence this common value
is finite, 
every $y \in \mathcal S_3$ satisfies
\begin{eqnarray*}
\Pi(\phi(\cdot) |y|^{\lambda}\,{\rm sgn }\,y) = \lambda \Pi(y) + \Pi(\phi).
\end{eqnarray*}
Therefore, if we choose $\phi$ so that $\Pi(\phi) > -n$, then (ii) is satisfied. This is the case of the Thomas--Fermi
equation; $\Pi(t^{-1/2}) = - 1/2 > -2$.\\
\\
\indent 
Under the assumptions we first show that a solution of Eq.(2) cannot belong to $\mathcal S_1$ and 
then we seek the possibility of solutions among functions in $\mathcal S_2$ or $\mathcal S_3$.\\
\\
{\bf Lemma 2.} {\it Under Assumption {\rm (i)} any function in $\mathcal S_1$ cannot be a solution of {\rm Eq.(2).}}\\
\\
{\it Proof.} The proof is by contradiction.
Assume that $y \in \mathcal S_1$ with domain $[a, \infty)$ is a solution of Eq.(2); $y^{(n)}(t) = F(t, y(t))$. 
By Assumption (i) we have $F(t, y(t)) = (y(t))^{\nu}g(t)$ for some function $g \in \mathcal S_1$ depending on $y$. 
Since the product of two functions in $\mathcal S_1$
also belongs to $\mathcal S_1$, it follows that $F(\cdot, y(\cdot))$ and $y^{(n)}$ also belong to $\mathcal S_1$. 
By applying Theorem 1 (iii) several times
we see that $y^{(i)} \in \mathcal S_1$ for $i = 0, 1, \ldots, n$. For a $v_1 > 0$ put $f = y'/y^{1+v_1} = - 1/v_1 (1/y^{v_1})'$. By integrating both sides
from $a$ to $t$ yields $1/(y(t))^{v_1} = c - v_1 \int_a^t f(\tau)d\tau$, where $c$ is a constant. Since $y \in \mathcal S_1$,
we have $\lim_{t \to \infty} y(t) = \infty$ and $\int_a^{\infty} f(\tau)d\tau$ is finite. Therefore, 
there exists a subset $J_1$ of $[a, \infty)$ whose Lebesgue measure is finite and only on which $|y'/y^{1+v_1}| \ge 1$ holds. 
Repeating this argument, for any $v_i > 0$ $(1 \le i \le n)$, there exists 
a subset $J_i$ of $[a, \infty)$ whose Lebesgue measure is finite and 
only on which $|y^{(i)}/(y^{(i-1)})^{1+v_i}| \ge 1$ holds. \\
\indent
For the number $\nu > 1$, making appropriate choices of $v_1, \ldots, v_n > 0$ and ${\mu}_1, \ldots, {\mu}_{n-1}> 1$, we see that
$y^{(n)}/y^{\nu}$ can be written as
\begin{eqnarray*}
\Big|\frac{y^{(n)}}{y^{\nu}}\Big| = \Big |\frac{y'}{y^{1 + v_1}} \Big|^{\mu_1}
\Big|\frac{y''}{{(y')}^{1 + v_2}} \Big|^{\mu_2} \cdots
\Big|\frac{y^{(n-1)}}{{(y^{(n-2)})}^{1 + v_{n-1}}} \Big|^{\mu_{n-1}}
\Big|\frac{y^{(n)}}{{(y^{(n-1)})}^{1 + v_n}}\Big|. 
\end{eqnarray*}
Using these numbers, we have for all $t \in [a, \infty) \backslash  J$, 
\begin{eqnarray*}
\Big|\frac{y^{(n)}(t)}{(y(t))^{\nu}}\Big| < 1,
\end{eqnarray*}
where $J = J_1 \cup J_2 \cdots \cup J_n$ and the Lebesgue measure of $J$ is finite. 
However, rewriting Eq.(2) as
\begin{eqnarray*}
\frac{y^{(n)}(t)}{(y(t))^{\nu}} = \frac{F(t,y(t))}{(y(t))^{\nu}},
\end{eqnarray*}
we see that this equality contradicts with the fact that the right-hand side of this equality is a function belonging to 
$\mathcal S_1$ by Assumption (i).
Hence any function in $\mathcal S_1$ cannot be a solution of Eq.(2). This completes the proof. \\
\\
\indent
In order to deduce more results on asymptotic solutions of Eq.(2), we employ Assumption (ii) on $F$.
They generalize the assertions of [1, Th8] in some aspects. \\
\\
{\bf Theorem 3.} {\it Let $y$ be a solution of {\rm Eq.(2)}. Then there exists an integer $i$ $(0 \le i \le n)$ such that 
for each integer $j$ $(i \le j \le n)$ an increasing sequence $\{t_k\}$ can be chosen such that 
$t_k \to \infty$ and $\lim_{k \to \infty} y^{(j)}(t_k) = 0$. 
}\\
\\
{\it Proof.} Recall that Theorem 1 tells us the following; a continuously differentiable function $y \in \mathcal S_2$ 
satisfies either $y' \in \mathcal S_2$ or $y' \in \mathcal S_3$ with $\Pi (y') = - \infty$, and moreover
a continuously differentiable function $y \in \mathcal S_3$ satisfies either
$y' \in \mathcal S_2$ or $y' \in \mathcal S_3$ with $\Pi (y') \le \Pi(y) - 1$.
Lemma 2 implies that a solution $y$ of Eq.(2) either belongs to $\mathcal S_2$ or to $\mathcal S_3$.
First consider the case $y \in \mathcal S_2$. Then we see that $y' \in \mathcal S_2$ or
$y' \in \mathcal S_3$ with $\Pi (y') = - \infty$. Repeating this argument we have, for $i = 0, 1, \ldots, n$, either
\begin{eqnarray*}
y^{(i)} \in \mathcal S_2 ~~{\rm or}~~ y^{(i)} \in \mathcal S_3 ~~{\rm with}~\Pi (y^{(i)}) = - \infty.
\end{eqnarray*}
Both mean that there exists an increasing sequence $\{t_k\}$ such that $t_k \to \infty$ and
$\lim_{k \to \infty} y^{(i)}(t_k) = 0$ for each $i$ $(0 \le i \le n)$. \\
\indent
Next assume $y \in \mathcal S_3$. Then we have either $y' \in \mathcal S_2$ or
$y' \in \mathcal S_3$ with $\Pi(y') \le \Pi(y) - 1$. Applying Theorem 1 again, in case of $y' \in \mathcal S_2$ it follows that
either $y'' \in \mathcal S_2$ or $y'' \in \mathcal S_3$ with $\Pi(y'') = - \infty$.  
In case of $y' \in \mathcal S_3$
we have either $y'' \in \mathcal S_2$ or $y'' \in \mathcal S_3$ with $\Pi(y'') \le \Pi(y') - 1 \le \Pi(y) - 2$. 
Continuing this process until $y^{(i)}$, we have the following two cases:
\begin{itemize}
\item[(1)]  either $y^{(i)} \in \mathcal S_2$ or $y^{(i)} \in \mathcal S_3$ with $\Pi(y^{(i)}) = - \infty$,
for some $i = 1, \ldots, n$; 
\item[(2)] $y^{(i)} \in \mathcal S_3$ for all $i = 0, 1, \ldots, n$ and 
$\Pi(y^{(n)}) \le \Pi(y) - n$.
\end{itemize}
When case (1) happens, the required assertion follows as above, particularly for all $j$ $(i \le j \le n)$.
When case (2) happens, due to Assumption (ii) and $c > - n$, we have
\begin{eqnarray*}
\mu \Pi(y) + c \le \Pi(F(\cdot, y(\cdot))) = \Pi(y^{(n)}) \le \Pi(y) - n,
\end{eqnarray*}
and 
\begin{eqnarray*}
    \Pi(y) \le - \frac{n + c}{\mu - 1} < 0.
\end{eqnarray*}
Moreover, for each $i = 1, \ldots, n$, it follows that
\begin{eqnarray*}
    \Pi(y^{(i)}) \le \Pi(y) - i \le - \frac{n + c}{\mu - 1}  - i < 0.
\end{eqnarray*}
Hence in case (2) we see that, for each $i$ $(0 \le i \le n)$, there exists an increasing sequence
for $y^{(i)}$ that satisfies the assertion of this theorem. This completes the proof. \\
\\
\indent
Let us return to Eq.(3) with appropriate functions $\phi$ such as both
$\Pi(\phi)$ and $\Xi(\phi)$ are finite. Suppose its solution $y$ is non-oscillatory.  Then we can
conclude that eventually $y^{(i)}(t)$ $(0 \le i \le n-1)$ tends to zero monotonically. Otherwise $y^{(i+1)}$ 
and hence $y^{(n)}$ change signs infinitely many times.
Due to the form of Eq.(3), so does $y$ itself, meaning that it becomes eventually oscillatory.
Therefore, any solution of Eq.(3) is eventually either oscillatory or tends to zero monotonically.\\
\\
%%%%%%%%%%%%%%%%%%%%%%%%%%%%%%%%%%%%%%%%%%%%%%%%%%%%%%%%%%%%%%%
{\bf References}
\begin{itemize}
\item[1.]  S. Tanimoto, Asymptotic behavior of functions and solutions of some nonlinear differential equations,
{\it Journal of Mathematical Analysis and Applications}, {\bf 138}, 511--521, 1989.
\item[2]  J. S. W. Wong, On the generalized Emden--Fowler equation, {\it SIAM Review}, {\bf 17}, 339--360,
1975.
\end{itemize}
\end{document}